\documentstyle[amscd]{amsart}
\numberwithin{equation}{section}
\theoremstyle{plain}
\newtheorem{theorem}{Theorem}[section]
\newtheorem{lemma}[theorem]{Lemma}
\newtheorem{cor}[theorem]{Corollary}
\newtheorem{prop}[theorem]{Proposition}
\newtheorem{claim}[theorem]{Conjecture}

\theoremstyle{definition}
\newtheorem{definition}[theorem]{Definition}
\newtheorem{example}[theorem]{Example}

\newtheorem{exmp}[theorem]{Problem}

\theoremstyle{remark}
\newtheorem{remark}[theorem]{Remark}

\numberwithin{equation}{section}

\begin{document}

\title{Motives and algebraic de Rham cohomology}

\author{Masanori Asakura}

\address{
Research Institute for Mathematics Sciences,
Kyoto University,
Oiwakecho, Sakyo-ku,
Kyoto, 606-8502, JAPAN
}
\email{asakura@@kurims.kyoto-u.ac.jp}

\thanks{
The author is supported by JSPS Research Fellowships for Young Scientists.
}

\subjclass{?}
\date{?}

\begin{abstract}
In this paper, we define a certain Hodge-theoretic structure 
for an arbitrary variety $X$ over the complex number field by using 
the theory of mixed Hodge module due to Morihiko~Saito. 
We call it an arithmetic Hodge structure of $X$. 
It is shown that 
extension groups of arithmetic Hodge structure do not vanish 
even for degree $\geq2$. 
Moreover, 
we define higher Abel-Jacobi maps from Bloch's higher Chow groups of 
$X$ to these extension groups. These maps
essentially involve the classical Abel-Jacobi maps by Weil and Griffiths,
and Mumford's infinitesimal invariants of 0-cycles on surfaces. 
\end{abstract}

\maketitle

%
\def\Spec{{\mathrm{Spec}}}
\def\Pic{{\mathrm{Pic}}}
\def\Alb{{\mathrm{Alb}}}
\def\Ext{{\mathrm{Ext}}}
\def\NS{{\mathrm{NS}}}
\def\Picv{{\mathrm{Pic^0}}}
\def\Div{{\mathrm{Div}}}
\def\CH{{\mathrm{CH}}}
\def\deg{{\mathrm{deg}}}
\def\rank{{\mathrm{rank}}}
\def\dim{{\mathrm{dim}}}
\def\codim{{\mathrm{codim}}}
\def\Coker{{\mathrm{Coker}}}
\def\ker{{\mathrm{ker}}}
\def\Image{{\mathrm{Image}}}
\def\Aut{{\mathrm{Aut}}}
\def\Hom{{\mathrm{Hom}}}
\def\Proj{{\mathrm{Proj}}}
\def\Sym{{\mathrm{Sym}}}
\def\Image{{\mathrm{Image}}}
\def\Gal{{\mathrm{Gal}}}
\def\GL{{\mathrm{GL}}}
\def\P{{\bold P}}
\def\C{{\bold C}}
\def\Q{{\bold Q}}
\def\Z{{\bold Z}}
\def\nb{\nabla}
\def\Lam{\Lambda}
\def\lam{\lambda}
\def\Om{\Omega}
\def\k{\kappa}
\def\l{\ell}
\def\z{\zeta}
\def\Te{\Theta}
\def\te{\theta}
\def\sg{\sigma}
\def\ve{\varepsilon}
\def\lra{\longrightarrow}
\def\ra{\rightarrow}
\def\la{\leftarrow}
\def\lla{\longleftarrow}
\def\Lra{\Longrightarrow}
\def\Lla{\Longleftarrow}
\def\da{\downarrow}
\def\hra{\hookrightarrow}
\def\lmt{\longmapsto}
\def\ot{\otimes}
\def\op{\oplus}
\def\wt#1{\widetilde{#1}}
\def\ol#1{\overline{#1}}
\def\us#1#2{\underset{#1}{#2}}
\def\os#1#2{\overset{#1}{#2}}
\def\lim#1{\us{#1}{\varinjlim}}
%
\def\et{\'{e}tale }
\def\vg{\varGamma}
\def\HS#1{{\mathrm{HS}}(#1)}
\def\QHS#1{{\mathrm{HS}}(#1)_{\Q}}
\def\MHS{{\mathrm{MHS}}}
\def\QMHS{\underline{{\mathrm{MHS}}}_{\Q}}
\def\MHM{{\mathrm{MHM}}}
\def\Perv{{\mathrm{Perv}}}
\def\VMHS{{\mathrm{VMHS}}}
\def\MM{\underline{{\mathrm{M}}}}
\def\sym{{\frak S}}
\def\O{{\mathcal O}}
\def\X{\Xi}

\def\Gr{{\mathrm{Gr}}}
\def\dego{\deg=0}
\def\Qb{\bar{\Q}}
\def\prim{0}

\def\al{\alpha}
\section{Introduction}
Let $X$ be a nonsingular projective surface over $\C$. 
D. Mumford was the first to show that the Chow group 
$\CH_0(X)$ of $0$-cycles 
on $X$ is more complicated than what was earlier believed (\cite{mumford}); 
in fact it can be ``enormous''. This has led to a difficulty
in understanding the precise nature of structure of 
the Chow group of zero-cycles. 
His theorem asserts that if the geometric genus of $X$ is not zero,
then the kernel $T(X)$ of the Albanese map
$$
A_0(X) \lra {\mathrm{Alb}}(X)
$$
cannot be finite-dimensional. Namely
$A_0(X)$ cannot  be given a 
geometric structure such as a complex torus.
Here $A_0(X)$ denotes the subgroup of $0$-cycles of degree $0$.

S. Bloch \cite{bloch} provided some
insight in this situation by conjecturing
that $T(X)$ is controlled by the transcendental part
of the cohomology $H^2(X,\Q)_{tr}=H^2(X,\Q)/\NS (X)$.
This is now called the Bloch 
conjecture (see Conjecture~\ref{Bloch} in \S\ref{anobservationtocycles}).
Combining this with the mixed motives as
introduced by P. Deligne and A. Beilinson,
his conjecture led to a conjectured existence of a 
filtration $F_{\mathcal M}$ on all of the higher Chow groups
$\CH^r(X,m;\Q) = \CH^r(X,m)\ot\Q$ (\cite{hc}), called the motivic
filtration. 
(In this paper we will simply write $\CH^r(X,m)$ instead
of $\CH^r(X,m;\Q)$, to be understood as Bloch's higher
Chow groups tensored with $\Q$.) This is fortified by
the following beautiful conjectural formula for any nonsingular 
projective variety $X$:
\begin{equation}\label{blochbeilinson}
\Gr_{F_{{\mathcal M}}}^{\nu}\CH^r(X,m)
=
\Ext_{{\mathcal M}}^{\nu}(\Q(0),H^{2r-m-\nu}(X)(r)).
\end{equation}
Here ${\mathcal M}$ is the category of mixed motives over 
$\Spec(\C)$ and $\Q(0)$ is the trivial motive.
It is proved that the motivic filtration is determined uniquely
if it exists (cf. \cite{jan}, \cite{ss}).
For example, it is conjectured that 
$F_{{\mathcal M}}^1\CH^r(X)=\CH^r(X)_{\hom}$ (the subgroup of homologically
trivial cycles), 
$F_{{\mathcal M}}^2\CH_0(X) = T(X)$, 
and $F_{{\mathcal M}}^{r+1}\CH^r(X)=0$.
In spite of much effort by many people, nobody has found a suitable
definition
of ${\mathcal M}$ and 
$F_{{\mathcal M}}^{\bullet}$ (see \cite{jan}, \cite{murre},
\cite{shuji}).

On the other hand, the category of the mixed motives ${\mathcal M}$ is
considered
to possess an exact faithful functor
${\mathcal M} \ra \MHS$ (called the {\it realization functor}) to
the category of graded polarizable $\Q$-mixed Hodge structures.
Therefore there should be the map
$$
\Gr_{F_{{\mathcal M}}}^{\nu}\CH^r(X)
\lra
\Ext_{\MHS}^{\nu}(\Q(0),H^{2r-\nu}(X)(r)).
$$
Unfortunately, the higher extension group $\Ext^{\nu}_{\MHS}$ 
($\nu \geq2$) always vanishes in the category $\MHS$.
This means that $T(X)=\Gr_{F_{{\mathcal M}}}^{2}\CH_0(X)$ for a surface $X$
cannot be captured by the extension of mixed Hodge structures.

\medskip

The purpose of this paper is to define a certain Hodge-theoretic structure 
for an arbitrary variety $X$ over the complex number field.
We call it an arithmetic Hodge structure of $X$. 
We define the decreasing filtration $F^{\bullet}\CH^r(X,m)$
on higher Chow groups and the
following map (which we call the {\it higher Abel-Jacobi map}):
$$
\rho^{\nu}_X:\Gr^{\nu}_F\CH^r(X,m)\lra 
\Ext^{\nu}_{\MM(\C)}(\Q(0),H^{2r-m-\nu}(X)(r)).
$$
Here $\MM(\C)$ denotes the category of arithmetic Hodge structure which
possesses the realization functor $r:\MM(\C)\ra \MHS$. 
When $\nu=1$, the map $\rho_X^1$ and 
the functor $r$ induces the usual Abel-Jacobi maps (cf. \S2.5).
When $\nu\geq2$, the higher Abel-Jacobi map $\rho^{\nu}_X$ gives a new tool
for Hodge theoretic study of algebraic cycles.
For example, it defines the {\it second Albanese map}
$$
T(X)\lra \Ext^2_{\MM(\C)}(\Q(0),H^{2n-2}(X)(n)),
$$
for a projective nonsingular variety $X$ of dimension $n$.
If it is injective, the Bloch conjecture is true (Theorem \ref{candq531}).
More generally,
we conjecture that $F^{\bullet}$ is discrete, that is, $F^N=0$ for $N\gg0$.
If our conjecture holds, 
our filtration $F^{\bullet}$ gives the motivic filtration
(cf. Theorem \ref{intro58} below).

\medskip

Let us explain arithmetic Hodge structures.
Let $X$ be a quasi-projective nonsingular variety over $\C$.
Then $X$ is defined by finitely many equations which have
finitely many coefficients.
By considering the coefficients as parameters of a space $S$,
we obtain the smooth family $f:X_S\ra S$.
Then we get the mixed Hodge module $H^{\bullet}(X_S/S)$ (cf. \cite{msaito}).
The usual mixed Hodge structure
$H^{\bullet}(X,\Q)$ appears on the fiber over a point of $S$.
The {\it arithmetic Hodge structure} is
defined to be the inductive limit of the
above mixed Hodge modules, $S$ running over the embedding
of the function field to $\C$.
The most important ingredient of arithmetic Hodge structure is the arithmetic
Gauss Manin connection on algebraic de Rham cohomology (see \S3.2).
Because of this connection,
the higher extension group does not necessarily vanish. 

\medskip

Our main result is as follows:
\begin{theorem}\label{intro58}
Let $X$ be a nonsingular projective variety over $\C$, and
$\CH^r(X,m)$ 
be the higher Chow group tensored with $\Q$. 
There is a decreasing filtration
$F^{\bullet}$ on $\CH^r(X,m)$ and a natural map
\begin{equation}\label{higherabel57}
\rho^{\nu}_X:\Gr^{\nu}_F\CH^r(X,m)\lra 
\Ext^{\nu}_{\MM(\C)}(\Q(0),H^{2r-m-\nu}(X)(r)),
\end{equation}
which has the following properties:
\begin{enumerate}
\renewcommand{\labelenumi}{(\theenumi)}
\item\label{fil531}
If $F^{r+1}\CH^r(X,m)=0$ for arbitrary 
$X$ and $r$, then $F^{\bullet}$ coincides with
the motivic filtration $F^{\bullet}_{\mathcal M}$. 
In particular, the vanishing implies
the Bloch conjecture $\ref{Bloch}$.
\item\label{filt531}
\begin{enumerate}
\item
$F^0\CH^r(X,m)=\CH^r(X,m)$. $F^1\CH^r(X)=\CH^r(X)_{\mathrm{hom}}$.
\item
$F^0\CH^r(X,m)=F^1\CH^r(X,m)$ for $m\geq1$.
\item
$F^2\CH^r(X,m)$ is contained in the kernel of the cycle map to
the Deligne cohomology group. These
coincide if the realization functor $r:\MM(\C)\ra \MHS$ is
fully faithful.
\item
$F^2\CH_0(X)=T(X)$.
\item
$F^{r+1}\CH^r(X,m)=F^{r+2}\CH^r(X,m)=\cdots$.
\end{enumerate}
\item\label{xiandde531}
There are the following natural maps:
$$
\Ext^{p}_{\MM(\C)}(\Q(0),H^{q}(X)(r))
\lra
\X_X^{r-p,q-r+p}(p)
\lra
\Lambda_X^{r-p,q-r+p}(p),
$$
where $\X_X^{p,q}(r)$ $($resp. $\Lambda_X^{p,q}(r)$ $)$ is 
defined as the cohomology of the following complex
$($see Definition $\ref{sh526}$ and Definition $\ref{mumdef58})$:
$$
F^{p+1}H^{p+q}_{dR}(X/\C)\ot \Omega^{r-1}_{\C/\Qb}
\os{{\nb}}{\ra}
F^{p}H^{p+q}_{dR}(X/\C)
\ot \Omega^r_{\C/\Qb}\\
\os{{\nb}}{\ra}
F^{p-1}H^{p+q}_{dR}(X/\C)
\ot \Omega^{r+1}_{\C/\Qb},
$$
$$(\text{resp. }
H^{q-1}(\Omega^{p+1}_{X/\C})\ot \Omega^{r-1}_{\C/\Qb}
\os{\ol{\nb}}{\ra}
H^{q}(\Omega^{p}_{X/\C})
\ot \Omega^r_{\C/\Qb}\\
\os{\ol{\nb}}{\ra}
H^{q+1}(\Omega^{p-1}_{X/\C})
\ot \Omega^{r+1}_{\C/\Qb}.
$$
In particular, the map \eqref{higherabel57} induces homomorphisms 
$$
\xi^{\nu}_X:\Gr^{\nu}_F\CH^r(X,m)\lra \X^{r-\nu,r-m}_X(\nu)
$$
and
$$
\delta^{\nu}_X:\Gr^{\nu}_F\CH^r(X,m)\lra \Lambda^{r-\nu,r-m}_X(\nu).
$$
In general, these maps have non-zero images, even when $\nu\geq2$, 
$($which follows, for example, from
Theorem $\ref{shujisaito})$.
\item
Assume that $X$ is defined over a number field $k$, and $X=X_k\ot_k\C$.
If $F^{\dim X+1}\CH_0(X)=0$, then $F^2\CH_0(X_k)=0$ and 
the rank of $\CH_0(X_k)$ is finite.
\end{enumerate}
\end{theorem}
We call the above map \eqref{higherabel57} the {\it $\nu$-th
higher Abel-Jacobi map}. Also, we remark that philosophically
speaking, the fact that $F^0\CH^r(X,m) = F^1\CH^r(X,m)$ for $m\geq1$
should not come as a surprise. Indeed, for $m\geq 1$ and by a 
standard weight argument, the cycle
class map from $\CH^r(X,m)$ to ordinary Betti cohomology is
zero. Thus $\CH^r(X,m)$ is already the kernel of a ``regulator''
map.

\medskip

By way of acknowledgement, a similar idea of arithmetic Hodge structures
was obtained by M. Green and P. Griffiths (\cite{G2}).
They defined the arithmetic Hodge structures
independently of the author\footnote
{There is a slight difference 
between their definition and ours.
In their definition, the datum of the $\Q$-structure does not appear.
Therefore the category of their arithmetic Hodge structures does not form
an abelian category, but only an exact category.}.
The terminology ``{\it arithmetic Hodge structure}" was introduced by them. 
It is appropriate to adopt their terminology in this paper.
However, our motivation for the above idea comes from the
the work of  Shuji Saito; that is,
the higher normal functions (\cite{AS}).
He defined the relative filtration $F_S^{\bullet}\CH^{\bullet}(X)$ and 
likewise $F_S^{\bullet}H_D^{\bullet}(X, \Q(\bullet))$ for 
a smooth projective family
$X \ra S$, and defined higher normal functions
$$
\Gr^{\nu}_{F_S}\CH^r(X) \lra
\vg(S, DJ_{X/S}^{r,\nu}).
$$
It was surprising for us that the above map has non-trivial
image. 

\medskip
This paper is organized as follows.

In \S2, we review Morihiko Saito's theory of mixed Hodge modules.
Here we omit the definition, because
what we need is not the precise definition but rather the formalism.
In \S3, we introduce the notion of arithmetic Hodge structure.
Moreover we define the spaces of Mumford's infinitesimal invariants, and
construct the natural maps from extension groups in $\MM(\C)$ to those
(Proposition \ref{muminv58}).
Finally in \S4, we construct the higher Abel-Jacobi maps, and study 
algebraic cycles, in particular, the Bloch conjecture.
\par
\vspace{0.25cm}
\hspace{-6.0mm} {\bf Acknowledgment}\par
\vspace{0.15cm}

We are especially indebted to Shuji Saito for his many helpful
suggestions.  It is my great pleasure to have had the opportunity
to work together with him.

We thank Shin-ichi Mochizuki, who kindly provided a proof for
the existence of the curve which appears in Theorem~\ref{asa68},
and Morihiko Saito, who has read through this paper carefully and
has come up with many helpful comments.

We also thank Sampei Usui and Takeshi Saito
for fruitful discussions and constant encouragement.

Special thanks are due to the two referees who have read through
the earlier versions of this paper pointing out many inaccuracies
and suggesting numerous improvements.

At last but not least, we would like to thank Noriko Yui and James D.
Lewis for their help for improving English and mathematical
presentations putting the paper into a presentable form to the reader.
Without their generous help and patience, this paper has never reached
this final form.

It is great pleasure to thank the organizers of the Banff Conference
on ``Arithmetic and Geometry of Algebraic Cycles'' for the excellent
conference.

\par
\vspace{0.25cm}
\hspace{-6.0mm} {\bf Notation and Conventions}\par
\vspace{0.15cm}
\begin{enumerate}
\item
A {\it variety} means a quasi-projective algebraic variety over a field.
We mainly work with algebraically closed fields of characteristic 0
(e.g. $\C$, $\Qb$).
\item
For a variety $X$ over $\C$, $X^{an}$ denotes the associated analytic
space: $X^{an}=X(\C)$.
\item
We always ignore torsion. We assume that all abelian groups
(e.g. $\CH^r(X,m)$ above) are tensored with $\Q$.
\item
In this paper, we fix an embedding $\Qb\hra\C$.
\end{enumerate}


\section{Morihiko Saito's 
formalism of mixed Hodge modules}\label{Hodgestructures}
In this section, we review the theory of mixed Hodge modules.
\subsection{Perverse sheaf}

A perverse sheaf is not a ``sheaf " in the usual sense, 
but rather a complex of sheaves with algebraically constructible
cohomology. 
In many ways it behaves like a sheaf (\cite{BBD}).

\medskip

Let $D^b_c(\Q_X)$ be the derived category of bounded complexes of 
$\Q_{X^{an}}$-sheaves 
with algebraically constructible cohomology for 
a complex algebraic variety $X$.
The {\it perverse $t$-structure} $(D^{\geq0},D^{\leq0})$ of $D^b_c(\Q_X)$
is the pair of additive full subcategories defined as follows:
\begin{enumerate}
\renewcommand{\labelenumi}{(\theenumi)}
\item
$
 K_{\Q}^{\bullet}\in D^{\leq0} 
\quad \Longleftrightarrow\quad
\dim {\mathrm{Supp}}H^i(K_{\Q}^{\bullet})\leq -i
\quad (\forall i),$
\item
$ K_{\Q}^{\bullet}\in D^{\geq0} 
\quad \Longleftrightarrow\quad
{\Bbb D}(K_{\Q}^{\bullet})\in D^{\leq0}
\quad (\forall i) $,
where $\Bbb D$ is the Verdier dual functor on $D^b_c(\Q_X)$.
\end{enumerate}
Here we define the dimension of empty set to be $-\infty$.
We call an object of $D^{\geq0}\cap D^{\leq0}$ a {\it perverse sheaf},
and denote its category by $\Perv(\Q_X)$.

The category of perverse sheaves is abelian.
There are the truncation functors ${}^p\tau_{\leq0}:D^b_c(\Q_X)\ra D^{\leq0}$
and ${}^p\tau_{\geq0}:D^b_c(\Q_X)\ra D^{\geq0}$, which are characterized as
the right and left adjoint of $D^{\leq0}\hra D^b_c(\Q_X)$ and
$D^{\geq0}\hra D^b_c(\Q_X)$ respectively.
We put ${}^pH^0:={}^p\tau_{\leq0}\circ {}^p\tau_{\geq0}(\simeq
{}^p\tau_{\geq0}\circ {}^p\tau_{\leq0})$, and call it
the perverse cohomology
functor.

There are the standard functors on $D^b_c(\Q_X)$:
\begin{equation}\label{gro56}
f_*,~f_!,~f^*,~f^!,~{\Bbb D},~\psi_g,~\phi_{g},~\ot,~\underline{\Hom}.
\end{equation}
We write ${}^pH^0f_*(-)[n]$ by ${}^pR^nf_*$ and so on.

The simplest example of perverse sheaf is a local system $L[\dim X]$
shifted by $\dim X$ when $X$ is nonsingular.
The category of local systems is not appropriate for the treatment of
homological algebra. 
Indeed, it is not closed under the above standard functors.
The category of constructible sheaves is closed under those,
but the one of perverse sheaves has more advantages.
Those properties of perverse sheaves are essentially used
to construct the theory of mixed Hodge modules.

\subsection{Regular holonomic $D$-modules}
Next we recall holonomic $D$-modules\footnote{
We use left $D$-modules, though M.Saito 
prefers right $D$-modules to it.} 
and regular singularities.

\medskip

Let $X$ be a nonsingular variety over a field $k$ of characteristic zero.

Let $M$ be an algebraic $D_X$-module which is quasi-coherent over $\O_X$. 
We denote $F_{i}D_X\subset D_X$ the differential order filtration, 
that is, it
is generated by differential operators 
$\sum_{\vert I\vert \leq i} f_I(x) \partial^I$.
An increasing filtration $F_{\bullet}M$ of coherent $\O_X$-submodules
is called a {\it good filtration} if 
\begin{enumerate}
\item
$\cup_i F_iM=M$, $F_{-k}M=0$ for some
integer $k\gg 0$, 
\item
$F_iD_X\cdot F_jM\subset F_{i+j}M$ for $\forall i,j$, 
\item
${\mathrm{gr}}^FM:=
\op_i{\mathrm{gr}}^F_iM$ is coherent over ${\mathrm{gr}}^FD_X$.
\end{enumerate}
Any coherent $D_X$-module has at least one good filtration.
Note that ${\mathrm{gr}}^FD_X$ is isomorphic to $\O_{T^*X}$
the structure sheaf of the
cotangent bundle $T^*X$ over $X$. The support of ${\mathrm{gr}}^FM$ in $T^*X$
is called the {\it characteristic variety} of $M$, which does not depend on
the choice of a good filtration.

In general, the dimension of the characteristic variety of a coherent
$D_X$-module is larger than or equal to $\dim X$ 
(Bernstein's inequality).
If the equality holds, we call the 
coherent $D_X$-module {\it holonomic}.

Let us recall the regularity of holonomic $D$-module\footnote
{There are several ways to define regularities of
holonomic $D$-modules. 
We follow the text of Z.Mebkhout \cite{meb} p.163 Proposition 5.4.2.}.
When $\dim X=1$, a holonomic $D_X$-module $M$ has {\it regular singularities}
if and only if there is a dense open set $U\subset X$ such that $M\vert_U$
is a connection with regular singularities, 
that is, there is a locally free sheaf $\ol{M}$
of finite rank over $\O_{\ol{U}}$ ($\ol{U}$ is the smooth completion of $U$) 
and log connection $\ol{\nb}:\ol{M}\ra \Omega_{\ol{U}/k}^1(\log D)\ot \ol{M}$
($D:=\ol{U}-U$) such that $(\ol{M},\ol{\nb})\vert_U\simeq (M\vert_U,\nb)$.
For $X$ of arbitrary dimension,
a holonomic $D_X$-module $M$ has {\it regular singularities} 
if and only if so does $H^0(i_C^*M)$ for
any $i_C:C\ra X$ with $C$ a nonsingular
curve over $k$, where $i_C^*$ is the derived inverse image functor
of $D$-modules.

\bigskip

The following theorem is due to Kashiwara and Mebkhout:
\begin{theorem}[Riemann-Hilbert correspondence]
Let $X$ be a nonsingular variety over $\C$.
Then the de Rham functor 
${\mathrm{DR}}:=\Omega_{X^{an}}^{\bullet}\ot(-)[\dim X]$ 
induces the equivalence
from the category of regular holonomic $D$-modules on $X$
to
the one of perverse sheaves with $\C$-coefficient:
$${\mathrm{DR}}:{\mathrm{Mod}}_{rh}(X)
\os{\sim}{\lra}\Perv(\C_X).$$
\label{kasmeb}
\end{theorem}

The above equivalence also induces the equivalence between local systems
and locally free $\O_X$-sheaves with integrable connection with
regular singularities (\cite{Deligne4}).

There are the standard operations as in $\eqref{gro56}$
on the derived category of bounded complexes of regular holonomic $D$-modules. 
Those operations are compatible under the Riemann-Hilbert correspondence.

\medskip

Again, let $X$ be a nonsingular variety over any field of characteristic zero.
${\mathrm{MF}}_{rh}(X)$ denotes the category of filtered
regular holonomic $D_X$-modules $(M,F)$, 
where $M$ is a regular holonomic $D_X$-module, and
$F$ is a good filtration.
Moreover
we denote by ${\mathrm{MFW}}_{rh}(X)$ the category of filtered
regular holonomic $D_X$-modules with weight filtration $(M,F,W)$,
where $W$ is a finite
increasing filtration of $D_X$-submodules (=the weight filtration).
These are not abelian, but exact categories such that the kernel
and cokernel objects exist for any morphisms.  
A complex in ${\mathrm{MF}}_{rh}(X)$ (resp. ${\mathrm{MFW}}_{rh}(X)$) 
$$(M',F)\lra (M,F)\lra (M'',F)$$
$$(\text{resp. }(M',F,W)\lra (M,F,W)\lra (M'',F,W))$$
is exact if and only if
$$
{\mathrm{gr}}^FM'\lra
{\mathrm{gr}}^FM\lra
{\mathrm{gr}}^FM''
$$
$$(\text{resp. }
{\mathrm{gr}}^F{\mathrm{gr}}^WM'\lra
{\mathrm{gr}}^F{\mathrm{gr}}^WM\lra
{\mathrm{gr}}^F{\mathrm{gr}}^WM''
)$$
is an exact sequence of sheaves of $\O_X$-modules.
\subsection{Mixed Hodge modules}
Morihoko Saito defined mixed Hodge modules, 
and proved the expected properties.

What is a mixed Hodge module?
Roughly speaking, it is a $\Q$-coefficient perverse sheaf with mixed Hodge
structure. Its notion contains {\it
 admissible variation of mixed Hodge structure}.
 
First we recall it:
\begin{definition}\label{advmhs530}
Let $X$ be a nonsingular variety over $\C$.
Then an admissible variation of mixed Hodge structures on $X$ 
is defined to be the data
$(H_{\Q},H_{\O},W_{\bullet}, F^{\bullet}, \nb, i)$ where:
\begin{itemize}
\item
$H_{\Q}$ is a local system of finite dimensional $\Q$-vector spaces on 
$X^{an}$,
\item
$H_{\O}$ is a 
locally free (Zariski) 
sheaf of ${\mathcal O}_X$-module of the same rank as $H_{\Q}$,
\item
$W_{\bullet}$ is a finite increasing filtration of $H_{\Q}$,
called the {\it weight filtration},
\item
$F^{\bullet}$ 
is a finite decreasing filtration on $H_{\O}$ 
by locally free ${\mathcal O}_X$-submodules,
called the {\it Hodge filtration},
\item
$\nb:H_{\O} \ra H_{\O}\us{{\mathcal O}_X}{\ot} \Om^1_{X/\C}$ 
an integrable connection
(called the {\it (algebraic) Gauss-Manin connection}),
\item
$i:H_{\Q}\ot\C \os{\sim}{\lra} \ker\nb^{an}$ 
(called the {\it comparison isomorphism}),
or equivalently, $i$ induces an isomorphism 
$H_{\Q}\ot\O_X^{an}\os{\sim}{\ra} H_{\O}^{an}$.
\end{itemize}
and these satisfy:
\begin{enumerate}
\renewcommand{\labelenumi}{(\theenumi)}
\item
For all points $s \in X^{an}$, the fiber 
$H_{\Q,s} \os{i}{\hookrightarrow}H_{\O} \ot \C(s)$ with the induced filtrations
$W_{\bullet,s}$ on $H_{\Q,s}$ and $F_s^{\bullet}$ on $H_{\O} \ot \C(s)$
defines a mixed Hodge structure. 
\item
(Griffiths transversality)
$W_{\bullet}$ and $F^{\bullet}$ satisfy the following:
$$
\nb(W_{\l}) \subset W_{\l}\ot \Om^1_{X/\C}, \quad
\nb(F^p) \subset F^{p-1}\ot \Om^1_{X/\C}\quad
\text{ for } \forall \l,p.
$$
\item
(polarizability)
For each $\l$, there are $\Q$-bilinear form 
$Q:\Gr^W_{\l}(H_{\Q}) \ot \Gr^W_{\l}(H_{\Q}) \ra \Q$
and ${\mathcal O}_X$-bilinear form 
$Q_X:\Gr^W_{\l}(H_{\O}) \ot\Gr^W_{\l}(H_{\O}) \ra {\mathcal O}_X$
satisfying:
\begin{enumerate}
\renewcommand{\labelenumi}{(\theenumi)}
\item
$Q$ and $Q_X$ are compatible under the comparison isomorphism $i$.
\item
$Q$ defines a polarization form on the $\Q$-Hodge structure
$(\Gr^W_{\l}(H_{\Q,s}),F^{\bullet}_s)$ for all $s \in X^{an}$.
\item
$Q_X(\nb(x),y)+Q_X(x,\nb(y))=dQ_X(x,y)$ for any local sections
$x,y \in H_{\O}$.
\end{enumerate}
\item(admissibility)
When the data
$(H_{\Q},H_{\O},W_{\bullet},F^{\bullet}, \nb, i)$ is
pulled back to a nonsingular complex algebraic curve $C$, it satisfies
the admissibility (\cite{kashiwara} 1.9):
\begin{enumerate}
\renewcommand{\labelenumi}{(\theenumi)}
\item
Any local monodromy around $\ol{C}-C$ is quasi-unipotent.
Here $\ol{C}$ denotes the smooth completion of $C$.
\item
The logarithm $N$ of the unipotent part of the local monodromy
admits a weight filtration relative to $W_{\bullet}$
(\cite{SZ} \S2).
\item
The Hodge filtration $F^{\bullet}$ can be extended to a locally
free subsheaf $\ol{F}^{\bullet}$ 
of the Deligne's canonical extension $\ol{H}_C$
such that $\Gr^{\bullet}_{\ol{F}}\Gr^{\ol{W}}_{\bullet}(\ol{H}_C)$ is
locally free.
\end{enumerate}
\end{enumerate}
We denote the category of admissible variations of mixed 
Hodge structures on $X$
by $\VMHS(X)$.
\end{definition}

In order to define mixed Hodge module, we replace $H_{\Q}$ by perverse sheaf, 
$(H_{\O},F^{\bullet},\nb)$ by filtered regular holonomic
$D_X$-module, and $i$ by the Riemann-Hilbert correspondence.

More precisely, 
\def\PervW{{\mathrm{PervW}}}
let $\PervW(\Q_X)$ be the category of perverse sheaves
$(K^{\bullet}_{\Q},W_{\Q})$ with finite increasing filtration $W_{\Q}$
of sub perverse sheaves
(=the weight filtration).
There is a natural functor from ${\mathrm{MFW}}_{rh}(X)$ to $\PervW(\C_X)$,
which maps $(M,F,W)$ to $({\mathrm{DR}}(M),{\mathrm{DR}}(W))$.
Denote the fiber product 
$\PervW(\Q_X)\times_{\PervW(\C_X)}{\mathrm{MFW}}_{rh}(X)$ by
${\mathrm{MFW}}_{rh}(X,\Q)$, that is, it consists
of the objects $(K_{\Q}^{\bullet},W_{\Q},M,F,W, i)$ where
$i:(K_{\Q}^{\bullet}\ot\C,W_{\Q}\ot\C)\simeq {\mathrm{DR}}(M,W)$ is
an isomorphism of filtered objects (cf. Theorem \ref{kasmeb}).

\begin{theorem}[Morihiko Saito \cite{msaito1}, \cite{msaito}]\label{SaM58}
There is an abelian full subcategory 
$$\MHM(X)\subset{\mathrm{MFW}}_{rh}(X,\Q)$$
for each nonsingular complex algebraic variety $X$, which satisfies the 
following properties.
We call an object of $\MHM(X)$ a mixed Hodge module on $X$.
$W$ $($resp. $F$ $)$ is called the weight filtration 
$($resp. Hodge filtration$)$.
\begin{enumerate}
\renewcommand{\labelenumi}{(\theenumi)}
\item 
The forgetful functor ${\mathrm{rat}}:\MHM(X)\ra \Perv(\Q_X)$ is
exact and faithful.
\item\label{strict526}
Any morphisms of mixed Hodge modules are strict with respect to Hodge
and weight filtration.
\item
Any pure weight full subcategory of $\MHM(X)$ $($called polarizable
Hodge modules$)$
is semi-simple. 
\item\label{grosam}
There are the standard operations on $D^b(\MHM(X))$ the derived category
of bounded complex of mixed Hodge modules:
$$f_*,~f_!,~f^*,~f^!,~{\Bbb D},~\psi_g,~\phi_{g,1},~\ot,~\underline{\Hom}.$$
Those functors 
satisfies the adjointness, projection formulas, and are compatible
with the ones on perverse sheaves under the forgetful functor ${\mathrm{rat}}$.
\item\label{adVMHS}
Let $\MHM(X)_{\mathrm{sm}}$ be the full subcategory of $\MHM(X)$ whose 
object consists of those whose
underlying perverse sheaf is smooth $($i.e. local system$)$.
Then there is the equivalence between the category of admissible
variations of mixed
Hodge structure $($
Definition $\ref{advmhs530})$ and $\MHM(X)_{\mathrm{sm}}$:
$$
\VMHS(X)\os{\sim}{\lra}
\MHM(X)_{\mathrm{sm}}.
$$
Here an admissible variation of mixed
Hodge structure 
$(H_{\Q},H_{\O},W_{\bullet}, F^{\bullet}, \nb, i)$ corresponds to
the mixed Hodge module $(H_{\Q}[\dim X], W_{\bullet-\dim X}, H_{\O}, F_{\bullet})$,
where $F_p:=F^{-p}$.
An admissible variation of Hodge structure of pure weight $n$ corresponds to
the Hodge module of weight $n+\dim X$.
In particular,
$\MHM(\Spec\C)$ is isomorphic to the category of 
graded polarizable $\Q$-mixed Hodge structures.
\item$($decomposition theorem$)$
Let $f:X\ra Y$ be a proper morphism between 
nonsingular varieties over $\C$.
Then there is a noncanonical isomorphism
$$f_*M\simeq \us{k}{\op}H^kf_*M[-k]
$$
in the derived category $D^b(\MHM(Y))$ for any pure weight Hodge module $M$.
\label{decom529}
\item
If $H$ is a mixed Hodge module of weight $\leq n$ $($resp. $\geq n$ $)$, 
then 
$H^kf_!(H)$ and $H^kf^*(H)$ 
$($resp. $H^kf_*(H)$, $H^kf^!(H)$, resp. ${\Bbb D}H$ $)$ 
are of weight $\leq k+n$
$($resp. $\geq k+n$, resp. $\leq -n$ $)$.
\end{enumerate}
\qed
\end{theorem}

Morihiko Saito firstly defined the category of polarizable Hodge modules
using vanishing cycle functor, 
and proved the stability of direct images of projective morphism,
Verdier dual and so on, which is the main result in \cite{msaito1}.
Next he defined mixed Hodge modules as a successive extension of
Hodge modules satisfying some stability conditions (\cite{msaito}).
The most important property is the existence of standard functors
(loc.cit. Theorem 4.3. etc). 
However, the precise construction of the category of mixed Hodge modules is 
long and complicated.
Moreover it is quite difficult to write it down explicitly.

\begin{definition}
$\Q_X(r)[\dim X]$ denotes the {\it Tate Hodge module} 
$$(\Q_X(r)[\dim X],W_{\bullet}, \O_X, F_{\bullet})$$
with $\Gr^W_{\dim X-2r}(\Q_X(r)[\dim X])=\Q_X(r)[\dim X]$ and
$\Gr^F_r(\O_X)=\O_X$,
(which is actually a mixed Hodge module of pure weight $\dim X-2r$
due to Theorem \ref{SaM58} \eqref{adVMHS}.)

We write the Tate Hodge module shifted by $-\dim X$ by $\Q_X(r)$:
$$
\Q_X(r):=(\Q_X(r)[\dim X])[-\dim X] \in D^b(\MHM(X)).
$$
\end{definition}
\subsection{Mixed Hodge modules over a field of characteristic 0}

Morihiko Saito has already pointed out 
that there are several modifications
of mixed Hodge modules.
Here we recall one of those. 

\begin{definition}\label{k61}
Let $k\subset\C$ be a subfield.
For a nonsingular variety $X_k$ over $k$, we define
$\MHM(X_k)$ to be the full subcategory of the fiber product
$${\mathrm{MFW}}_{rh}(X_k)\times_{{\mathrm{MFW}}_{rh}(X_{\C})}
\MHM(X_{\C})\qquad (X_{\C}:=X_k\ot\C)$$
such that the polarization on each weight quotient is defined over $k$.
\end{definition}

Since all operations such as \eqref{gro56}
can be also defined on ${\mathrm{MFW}}_{rh}(X_k)$, 
we have those induced on  
$\MHM(X_k)$. Each property of those operations can be reduced to
those of $\MHM(X_{\C})$.
In particular, Theorem \ref{SaM58} also holds, if
we replace $\MHM(X_{\C})$ by $\MHM(X_k)$, and
$\VMHS(X_{\C})$ by the full subcategory
of ${\mathrm{MFW}}_{rh}(X_k)_{\mathrm{sm}}
\times_{{\mathrm{MFW}}_{rh}(X_{\C})_{\mathrm{sm}}}
\VMHS(X_{\C})$ such that the polarization on each weight quotient
is defined over $k$.
(${\mathrm{MFW}}_{rh}(X_k)_{\mathrm{sm}}$ denotes the full subcategory
of ${\mathrm{MFW}}_{rh}(X_k)$ generated by objects whose underlying
$D$-module is a locally free $\O_X$-module of finite rank.)
For more details, see his exposition (\cite{formal} (1.8) (ii)).

\begin{cor}\label{leray58}
Let $f:X\ra Y$ be a proper morphism between nonsingular varieties over $k$.
Then there is the Leray spectral sequence
\begin{equation}\label{leray529}
E_2^{pq}=\Ext^p_{\MHM(Y)}(\Q_Y(0),H^qf_*\Q_X(r))\Longrightarrow
\Ext^{p+q}_{\MHM(X)}(\Q_X(0),\Q_X(r)),
\end{equation}
which degenerates at $E_2$ terms.
\end{cor}
\begin{pf}
The existence of the spectral sequence \eqref{leray529} follows from
the existence of the standard functors (Theorem \ref{SaM58} \eqref{grosam})
and the fact $f^*\Q_Y(r)=\Q_X(r)$.
The $E_2$-degeneration follows from 
the decomposition theorem (Theorem \ref{SaM58} \eqref{decom529}).
\end{pf}
\begin{cor}[\cite{formal} (8.3)]
Let $X$ be a nonsingular variety over $k$. Then
there is a cycle map from the higher Chow group to
the extension group of mixed Hodge modules:
\begin{equation}\label{hcycle529}
\CH^r(X,m)\lra \Ext^{2r-m}_{\MHM(X)}(\Q_X(0),\Q_X(r)).
\end{equation}
\label{cycle529}
\end{cor}
\begin{pf}
(Here we recall the construction of the cycle map \eqref{hcycle529} 
for the convenience of the reader.)

Let $Y$ be a nonsingular variety over $k$ of dimension $d_Y$, 
and $i_Z:Z\hra Y$ a
closed subscheme (not necessarily smooth nor irreducible) of pure
dimension $d_Z$.
Let $j_Z:U\hra Y$ be the complement of $Z$.
We define $i_{Z*}i_Z^!\Q_Y(r)=
{\mathrm{Cone}}(\Q_Y(r)\ra j_{Z*}j^*_Z\Q_Y(r))[-1]$,
and $i_{Z*}\Q_Z(r)={\mathrm{Cone}}(j_{Z!}j^*_Z\Q_Y(r)\ra \Q_Y(r))$,
which are complexes of mixed Hodge modules on $Y$ 
by Theorem \ref{SaM58} \eqref{grosam}.
Note that ${\mathrm{rat}}(i_{Z*}\Q_Z(r))=i_{Z*}\Q_{Z^{an}}$
and ${\mathrm{rat}}(i_{Z*}i^!_Z\Q_Y(r))=i_{Z*}i_Z^!\Q_{Y^{an}}$.

\medskip

Let $Z_k$ be an irreducible component of $Z$, and
$z=\sum_k n_k Z_k$ be a cycle on $Y$ with $n_k\in \Q$.
First we construct the natural morphism
\begin{equation}\label{step1}
c_z:i_{Z*}\Q_Z(0)\lra \Q_Y(d_Y-d_Z)[2d_Y-2d_Z] 
\end{equation}
associated to the cycle $z$ in the derived category $D^b(\MHM(Y))$
as follows.
Let $h_k:\wt{Z}_k\ra Z_k$ be a resolution of singularity, and
$h:\amalg \wt{Z}_k\ra Y$ the composition with the inclusion $i_Z$. 
There are the pull-back morphisms 
$\alpha:i_{Z*}\Q_Z(0)\ra i_{Z*}h_{k*}
\Q_{\wt{Z}_k}(0)^{\op}=h_*\Q_{\wt{Z}_k}(0)^{\op}$
and $\beta:\Q_Y(0)\ra h_*\Q_{\wt{Z}_k}(0)^{\op}$.
Applying the Verdier dual $\Bbb D$ on
the latter morphism $\beta$, we have 
$h_*\Q_{\wt{Z}_k}(d_Z)^{\op}[2d_Z]\ra \Q_{Y}(d_{Y})[2d_{Y}]$
because $Y$ and $Z_k$ are nonsingular.
Composing this with
$\op n_k: h_*\Q_{\wt{Z}_k}(0)^{\op}\ra  
h_*\Q_{\wt{Z}_k}(0)^{\op}$ and
the morphism $\alpha$, we have \eqref{step1}.
It is a routine work to show its 
independency of the choice of resolutions of $Z_k$.

Put $Y=X$, and compose the morphism \eqref{step1} with $\Q_X(0)\ra
i_{Z*}\Q_Z(0)$.
Then we have 
the cycle map \eqref{hcycle529} for $m=0$. 

\medskip

Next we construct the cycle map \eqref{hcycle529} for $m\geq1$. 
We write
$$
\Delta^m=\Spec k[t_0,\cdots, t_m]/(\sum_{i=0}^m t_i-1)
$$
and its faces by
$$
\Delta_I: t_{i_0}=\cdots=t_{i_k}=0
$$
for a subset $I=\{i_0,\cdots,i_k\}\subset \{0,\cdots,m\}$.
Put $Y=X\times \Delta^m$, and $i:X\times (\cup_i \Delta_i)
\hra Y$ the closed immersion and $j:X\times U\hra Y$ its
complement.
Let $z=\sum n_kZ_k$ be a cycle on $Y$ (where $Z_k$ are reduced
and irreducible subvarieties of dimension $d_Y-r$) which intersects with
all faces $X\times \Delta_I$ properly, and 
\begin{equation}\label{inter61}
z\cdot[X\times \Delta_i]=0
\quad \text{for }0\leq\forall i\leq m.
\end{equation}
Put $Z=\cup_kZ_k$ the support of $z$.

We first claim:
\begin{lemma}\label{step2}
The morphism $c_z$ \eqref{step1} is uniquely lifted to
the morphism
$$
i_{Z*}\Q_Z(0)\lra j_!j^*\Q_Y(r)[2r] .
$$
\end{lemma}
In fact, using the distinguished triangle
$$
j_!j^*\Q_Y(r)\lra \Q_Y(r)\lra 
i_*\Q_{X\times(\cup_i \Delta_i)}(r)\os{+1}{\lra},
$$
it follows from
\begin{equation}\label{step3}
\Ext^{\nu}_{\MHM(Y)}(i_{Z*}\Q_Z(0),i_*\Q_{X\times(\cup_i \Delta_i)}(r))
=0 \quad\text{for }\forall \nu<2r,
\end{equation}
and
\begin{equation}\label{step4}
\text{
the composition of $c_z$ with  $\Q_Y(r)\lra 
i_*\Q_{X\times(\cup_i \Delta_i)}(r)$ is zero.}
\end{equation}
Since $ i_*\Q_{X\times(\cup_i \Delta_i)}(r)$ is isomorphic to 
$$
\os{m}{\us{i=0}{\op}}\Q_{X\times \Delta_i}(r)\lra
\us{i<j}{\op}\Q_{X\times \Delta_{ij}}(r)\lra\cdots
\lra \us{\vert I\vert =m}{\op}\Q_{X\times \Delta_I}(r),
$$
\eqref{step3} and \eqref{step4} follow from
\begin{equation}\label{step5}
\Ext^{\nu}_{\MHM(Y)}(i_{Z*}\Q_Z(0),i_{I}*\Q_{X\times\Delta_I}(r))
=0 \quad\text{for }\forall I\text{ and }\forall \nu<2r
\end{equation}
where $i_I:X\times\Delta_I\hra Y$, and
\begin{equation}\label{step6}
\text{
the composition of $c_z$ with  $\Q_Y(r)\lra 
i_*\Q_{X\times\Delta_i}(r)$ is zero for $0\leq\forall i\leq m$.}
\end{equation}
\eqref{step6} follows from \eqref{inter61}. 
In order to see \eqref{step5}, we have
\begin{align*}
\Ext^{\nu}_{\MHM(Y)}(i_{Z*}\Q_Z(0),i_{I*}\Q_{X\times\Delta_I}(r))
&=\Ext^{\nu}_{\MHM(X\times \Delta_I)}
(i_{Z_I*}\Q_{Z_I}(0),\Q_{X\times\Delta_I}(r))\\
&=\Ext^{\nu}_{\MHM(X\times \Delta_I)}
(i_{Z_I*}\Q_{Z_I}(0),i_{Z_I*}i_{Z_I}^!\Q_{X\times\Delta_I}(r))
\end{align*}
by the adjunction,
where $Z_I:=Z\cap( X\times\Delta_I)$.
Then \eqref{step5} follows from the fact that
$H^k((i_{Z_I*}\Q_{Z_I}(0))=0$ for $k>d_{Z_I}$
and $H^k(i_{Z_I*}i_{Z_I}^!\Q_{X\times\Delta_I}(r))=0$
for  $k<d_{X\times \Delta_I}+r$, which can be checked
on the underlying perverse sheaves (cf. \cite{formal} (7.15)).
Thus we have shown Lemma \ref{step2}.
\begin{lemma}\label{step7}
Let $p:Y\ra X$ be the projection. Then
$ p_*j_!j^*\Q_Y(r)= \Q_X(r)[-m]$.
\end{lemma}
In fact, 
since $j=1\times j_0:X\times (\Delta^m\setminus \cup_i
\Delta_i)\hra Y=X\times \Delta^m$, we can assume that $X=\Spec k$.
Then the assertion follows from \cite{formal} (7.6).

\medskip

By Lemma \ref{step2} and Lemma \ref{step7}, we have the morphism
$$
\Q_X(0)=p_*\Q_Y(0)\ra p_* i_{Z*}\Q_Z(0)
\ra   p_*j_!j^*\Q_Y(r)[2r]= \Q_X(r)[2r-m],
$$
which gives the cycle map \eqref{hcycle529}.
\end{pf}

\subsection{The Carlson isomorphisms}\label{carlson}
The category $\MHS$ is not semi simple, and the Yoneda extension groups
in $\MHS$ have the well known explicit description, due to Carlson (
\cite{jcarlson}).

Here we review it.
Since $\MHS$ carries the neutral Tannakian structure, we have:
$$
\Ext_{\MHS}^{\bullet}(H_1,H_2) =
\Ext_{\MHS}^{\bullet}(\Q(0),H_1^*\ot H_2).
$$
\begin{theorem}[Carlson]
Let $H=(H_{\Q},W_{\bullet},F^{\bullet})$ be a graded polarizable 
$\Q$-mixed Hodge structure. 
Then
\begin{enumerate}
\renewcommand{\labelenumi}{(\theenumi)}
\item
$
\Ext_{\MHS}^{1}(\Q(0),H)=
\Ext_{\MHS}^{1}(\Q(0),W_0H)=W_{-1}H_{\C}/W_{-1}H_{\C}\cap
(F^0W_0H_{\C}+W_0H_{\Q})
$
Here an element $\xi \in W_{-1}H_{\C}$ corresponds to the following extension
of graded polarizable mixed Hodge structures:
$$
0 \lra W_0H \lra \wt{H} \lra \Q(0) \lra 0,
$$
where $\wt{H}=(\wt{H}_{\Q}, W_{\bullet}, F^{\bullet})$ is the mixed Hodge
structure with $\wt{H}_{\Q}=W_0H_{\Q}\op \Q(0)$, and the weight
filtration $W_0\wt{H}_{\Q}=\wt{H}_{\Q}$, $W_{\l}\wt{H}_{\Q}=W_{\l}H_{\Q}$
$(\l \leq -1)$, and the Hodge filtration 
$F^p\wt{H}_{\C}=F^pW_0H_{\C}$ $(p \geq 1)$, 
$F^q\wt{H}_{\C}=F^qW_0H_{\C}+F^0\wt{H}_{\C}$ $(q \leq -1)$ and 
$$
F^0\wt{H}_{\C}=F^0W_0H_{\C}+ \C \cdot (\xi \op 1).
$$
\item
$\Ext_{\MHS}^{\nu}(\Q(0),H)=0$ for $\nu \geq 2$.
$($This is a formal consequence of
the fact that the functor $\Ext^1(\Q(0), -)$ is right exact, 
which follows from
the above explicit description.$)$
\end{enumerate}
\label{carl}
\end{theorem} 
\begin{remark}
The author learned the above formula from Morihiko Saito, 
which is more complicated than the one by J. Carlson. 
Originally, J. Carlson calculated the Yoneda extension 
groups in the category
of mixed Hodge structures which are not necessarily graded polarizable, 
and obtained a simple description (\cite{jcarlson}),
$$
\Ext^1(\Q(0), H)=W_0H_{\C}/F^0W_0H_{\C}+W_0H_{\Q}.
$$
\end{remark}
\medskip

In particular, if $H=H^{2r-1}(X^{an},\Q(r))$ for a complex projective
nonsingular variety $X$, 
\begin{equation}\label{carliso}
\Ext_{\MHS}^{1}(\Q(0),H)=J^r(X)_{\Q},
\end{equation}
where $J^r(X)_{\Q}=H_{\C}/F^0H_{\C}+H_{\Q}$ is the {\it $r$-th intermediate
Jacobian} of $X$.
Using the mixed Hodge modules,
we can define the Abel-Jacobi map
\begin{equation}\label{AJ}
\rho:\CH^r(X)_{\hom} \lra J^r(X)_{\Q},
\end{equation}
where $\CH^r(X)_{hom}$ denotes the $\Q$-tensored Chow group of $X$
generated by algebraic cycles of codimension $r$ which are homologically
equivalent to $0$.

Let
\begin{equation}
\CH^r(X)\lra \Ext^{2r}_{\MHM(X)}(\Q_X(0),\Q_X(r)).
\end{equation}
be the cycle map as in Corollary \ref{cycle529}.
By the Leray spectral sequence for the structure morphism $f:X\ra \Spec\C$
(Corollary  \ref{leray58}) and 
the vanishing of $\Ext_{\MHS}$ of deree $\geq2$, 
we get the following commutative diagram:
\begin{equation}\label{algab33}
\begin{matrix}
0& &0\\
\downarrow && \downarrow\\
\CH^r(X)_{\hom}&\os{\rho}{\lra}&\Ext^1_{\MHS}(\Q(0),H^{2r-1}(X,\Q(r)))\\
\downarrow && \downarrow\\
\CH^r(X)&\lra&\Ext^{2r}_{\MHM(X)}(\Q_{X}(0),\Q_{X}(r))\\
\downarrow && \downarrow\\
\CH^r(X)/\CH^r(X)_{\hom}&\lra& H^{2r}(X,\Q)\cap H^{r,r}\\
\downarrow && \downarrow\\
0& &0.
\end{matrix}
\end{equation}
Then the top horizontal arrow $\rho$ is defined to be the Abel-Jacobi map.

\bigskip

The Abel-Jacobi maps were originally defined by A.Weil and P.Griffiths.
Let us recall their definition.

Let $Z\in \CH^r(X)_{\hom}$ be an algebraic cycle.
There is a topological $(2\dim X-2r+1)$-cycle $\Gamma$ whose boundary
is $Z$: $\partial \Gamma =Z$.
Then the Abel-Jacobi class is defined as follows: 
\begin{equation}\label{class33}
\rho(Z)=\sum_{j=1}^g\left(\int_{\Gamma}\omega_j\right)\omega_j^*\in J^r(X),
\end{equation}
where $\omega_1,\cdots,\omega_g\in F^{n-r+1}H^{2n-2r+1}(X,\C)$ is a basis,
and $\omega_1^*,\cdots,\omega_g^*\in H^{2r-1}(X,\C)/F^r$ 
denotes the Serre dual
class of those: $\langle \omega_i,\omega_j^*\rangle=\delta_{ij}$.

It is known that
the both definitions of Abel-Jacobi maps
coincide (cf. \cite{EZ}, \cite{jan3}).
\section{Arithmetic Hodge structure}
Let $X$ be a nonsingular algebraic variety over $\C$.
When we study the algebraic cycles of $X$ (in particular, contained in 
the kernel of the Abel-Jacobi map), the many difficulties are often caused
from the fact of the vanishing of extension groups of mixed Hodge structure
of degree $\geq2$.

In this section, we introduce the notion of arithmetic Hodge structure
(or more generally, arithmetic Hodge module), which has non-trivial higher
extension groups. 

\subsection{}
Let $X$ be a quasi-projective nonsingular variety over $\C$.
Then $X$ is defined by finitely many equations which possess
finitely many coefficients.
By considering the coefficients as parameters of a space $S$,
we can obtain a model $f_S:X_S\ra S$ and the Cartesian diagram:
\begin{equation}\label{model}
\begin{CD}
X_S @<<< X\\
@V{f_S}VV  @VVV\\
S @<{a}<<\Spec \C,
\end{CD}
\end{equation}
where $S$ is a nonsingular variety over $\Qb$, and the map $a$ factors
through the generic point
$\Spec~\Qb(S) \hookrightarrow S$.
\begin{definition}
We define the abelian category
$$\MM(X)=\lim{X_S} \MHM(X_S).$$
Here $\MHM(X_S)$ is the category of mixed Hodge modules on the variety $X_S$
over $\Qb$. In the above limit,
$X_S$ runs over all models \eqref{model}, and for a morphism
$j:X_{S'}\ra X_S$ of the models, we take the pull-back 
$j^*:\MHM(X_S)\ra\MHM(X_{S'})$.
We call it the category of {\it arithmetic Hodge modules}.
In particular, we call $\MM(\C):=\MM(\Spec\C)$
the category of {\it arithmetic Hodge structures}.
\end{definition}

$\MM(X)$ is an abelian category.
By Theorem \ref{SaM58} \eqref{grosam} and 
the remark about Definition \ref{k61}, 
the category of arithmetic Hodge modules
carries the standard functors:
$$f_*,~f_!,~f^*,~f^!,~{\Bbb D},~\ot,~\underline{\Hom}.$$

The {arithmetic Tate Hodge module} $\Q_X(r)$ is defined as the equivalence
class of $[\Q_{X_S}(r)]$ of the Tate Hodge module on a model $X_S$.
It does not depend on the choice of a model $X_S$. It follows from the fact
that for
any two models $X_S$ and $X_{S'}$, there is another model 
$X_{S^{\prime\prime}}$ with
morphisms $X_{S^{\prime\prime}}\ra X_S$ and 
$X_{S^{\prime\prime}}\ra X_{S'}$ of models.
We write $H^{\bullet}(X,\Q(r)):=R^{\bullet}f_*\Q_X(r)$ (or, simply
$H^{\bullet}(X)(r)$) for the structure morphism $f:X\ra\Spec\C$.

\subsection{The space of Mumford's infinitesimal invariants}
Yoneda extension groups in the category of arithmetic Hodge structures are
complicated. I do not know its explicit form in general.
Here we construct a natural map from the extension groups to the spaces of
Mumford's infinitesimal invariants, which have explicit forms.

\medskip

Let $X$ be a nonsingular projective variety over $\C$.
Let us recall the {\it arithmetic Gauss-Manin connection} on the algebraic
de Rham cohomology (cf.\cite{H}):
\begin{equation}\label{AGM}
\nb :H_{dR}^q(X/\C)\lra H_{dR}^q(X/\C)\ot_{\C}\Omega^1_{\C/\Qb}.
\end{equation}
Let $f:X\ra \Spec\C$ be the structure morphism.
The exact sequence
$$
0\lra f^*\Omega_{\C/\Qb}^1
\lra \Omega^{1}_{X/\Qb}
\lra \Omega^{1}_{X/\C}
\lra 0
$$
induces the following exact sequence
\begin{equation}
0\lra \Omega^{\bullet-1}_{X/\C}\ot f^*\Omega_{\C/\Qb}^1
\lra \Omega^{\bullet}_{X/\Qb}/U^2
\lra \Omega^{\bullet}_{X/\C}
\lra 0,
\label{le79}
\end{equation}
where $U^2:=\Image(\Omega^{\bullet-2}_{X/\Qb}\ot f^*\Omega_{\C/\Qb}^2\ra
\Omega^{\bullet}_{X/\Qb})$. Applying $Rf_*$ on \eqref{le79},
we obtain the coboundary map $H^q_{dR}(X/\C):=R^qf_*\Omega^{\bullet}_{X/\C}
\ra R^qf_*(\Omega^{\bullet-1}_{X/\C}\ot f^*\Omega_{\C/\Qb})\simeq
H^q_{dR}(X/\C)\ot \Omega^1_{\C/\Qb}$ (where the last isomorphism is the 
projection formula), 
and hence the arithmetic Gauss-Manin connection
\eqref{AGM}.

By definition, the arithmetic Gauss-Manin connection \eqref{AGM}
satisfies the Griffiths transversality
$$
\nb(F^p)\subset F^{p-1}\ot\Omega^1_{\C/\Qb},
$$
where $F^{s}=F^{s}H_{dR}^q(X/\C)=H^q(X,\Omega_{X/\C}^{\bullet\geq s})$ 
denotes the Hodge filtration.
\begin{definition}\label{sh526}
Let $p,q,r$ be integers.
Then $\X_X^{p,q}(r)$
is defined as the cohomology at the middle term of the following complex 
induced from the arithmetic Gauss-Manin connection \eqref{AGM}:
$$
F^{p+1}H^{p+q}_{dR}(X/\C)\ot \Omega^{r-1}_{\C/\Qb}
\os{{\nb}}{\ra}
F^{p}H^{p+q}_{dR}(X/\C)
\ot \Omega^r_{\C/\Qb}\\
\os{{\nb}}{\ra}
F^{p-1}H^{p+q}_{dR}(X/\C)
\ot \Omega^{r+1}_{\C/\Qb}.
$$
\end{definition}

We construct a natural map from the extension groups in $\MM(\C)$ to
the space $\X_X^{p,q}(r)$.

\medskip

Let $X_S\to S$ be the model of $X$ \eqref{model}.
Recall the category $\MHM(S)$ is defined to be a full subcategory of
the fiber product
${\mathrm{MFW}}_{rh}(S)
\times_{{\mathrm{MFW}}_{rh}(S_{\C})}
\MHM(S_{\C})$ such that the polarization on each weight quotient is 
defined over $\Qb$.
There is the forgetful functor $\MHM(S)\ra 
{\mathrm{MF}}_{rh}(S)$, which is exact by Theorem  \ref{SaM58} 
\eqref{strict526}:
\begin{equation}
\MHM(S)\lra {\mathrm{MF}}_{rh}(S),
\quad
(K^{\bullet}_{\Q},M,F,W)\mapsto (M,F).
\end{equation}
Therefore we have the well-defined map of Yoneda extension groups:
\begin{equation}\label{graded58}
\Ext_{\MHM(S)}^p(\Q_S(0),R^qf_*\Q_{X_S}(r))
\lra
\Ext_{{\mathrm{MF}}_{rh}(S)}^p
(\O_S,H^q_{dR}(X_S/S)(r))
\end{equation}
Note that ${\mathrm{MF}}_{rh}(S)$ is an exact category (see. \S2.2), 
and we can define its 
derived category $D^b({\mathrm{MF}}_{rh}(S))$ (cf. \cite{BBD} 1.1). 
Then we can describe
$\Ext_{{\mathrm{MF}}_{rh}(S)}^p(A,B)=
\Hom_{D^b({\mathrm{MF}}_{rh}(S))}(A,B[p])$.
We will take the inductive limit over $S$
of the right hand side in \eqref{graded58}.
So we may assume $S=\Spec\O$ where $\O$ is a regular local ring over $\Qb$
with regular parameter $x_1,\cdots,x_n$.
Then we can see $D_S=\O[\partial_1,\cdots,\partial_n]$ where $\partial_i$
is the derivation on $\O$ over $\Qb$ such that
$[\partial_i,x_j]=\delta_{ij}$.
Let $D_S(m)$ be the filtered $D$-module $(D_S,F_{\bullet+m})$
with the differential order filtration shifted by $m$.

\begin{lemma}\label{van527}
$\Ext_{{\mathrm{MF}}_{rh}(S)}^p
(D_S(m),M)=0
\quad \forall p\geq1.$
\end{lemma}
\begin{pf}Standard.
\end{pf}

There is the Koszul resolution of $\O_S$
\begin{equation}\label{kos527}
0\ra D_S(-n)\ot_{\Qb}\os{n}{\wedge}V\ra \cdots \ra D_S(-1)\ot_{\Qb}V
\ra D_S\ra \O_S\ra0,
\end{equation}
where $V$ is a $n$-dimensional $\Qb$-vector space with a basis 
$e_1,\cdots,e_n$,
and the boundary map
$D_S(-p)\ot_{\Qb}\os{p}{\wedge}V\ra D_S(-p+1)\ot_{\Qb}\os{p-1}{\wedge}V$
is given by $\xi\ot e_{i_1}\wedge\cdots \wedge e_{i_p}\mapsto
\sum_{j=1}^{p}(-1)^j\xi\cdot\partial_{i_j}\ot
e_{i_1}\wedge\cdots \wedge\hat{e}_{i_j}\wedge\cdots\wedge e_{i_p}$.
It is also the exact sequence in the exact category ${\mathrm{MF}}_{rh}(S)$.

Using the exact sequence \eqref{kos527} and applying Lemma \ref{van527},
we can see that the right hand side of \eqref{graded58} is isomorphic to
the cohomology at the middle term of the following complex 
induced from the algebraic Gauss-Manin connection:
$$
F^{r-p+1}H_{dR}^q(X_S/S)\ot \Omega^{p-1}_{S/\Qb}
\os{{\nb}}{\ra}
F^{r-p}H_{dR}^q(X_S/S)
\ot \Omega^p_{S/\Qb}
\os{{\nb}}{\ra}
F^{r-p-1}H_{dR}^q(X_S/S)
\ot \Omega^{p+1}_{S/\Qb}.
$$

Passing to the limit over $S$, we thus have:
\begin{prop}\label{xi527}
Let $X$ be a nonsingular projective variety over $\C$.
Then there is the following natural map:
$$
\Ext_{\MM(\C)}^p(\Q(0),H^q(X)(r))\lra \X_X^{r-p,q-r+p}(p).
$$
\end{prop}

\bigskip

\begin{definition}\label{mumdef58}
Let $p,q,r$ be integers.
The space of {\it Mumford's infinitesimal invariants} 
$\Lambda_X^{p,q}(r)$
is defined as the cohomology at the middle term of the following complex 
induced from the arithmetic Gauss-Manin connection \eqref{AGM}:
$$
H^{q-1}(\Omega^{p+1}_{X/\C})\ot \Omega^{r-1}_{\C/\Qb}
\os{\ol{\nb}}{\ra}
H^{q}(\Omega^{p}_{X/\C})
\ot \Omega^r_{\C/\Qb}\\
\os{\ol{\nb}}{\ra}
H^{q+1}(\Omega^{p-1}_{X/\C})
\ot \Omega^{r+1}_{\C/\Qb}.
$$
\end{definition}
Clearly, there is the natural map,
\begin{equation}\label{xtol}
\X_X^{p,q}(r)\lra \Lambda_X^{p,q}(r).
\end{equation}

By the similar argument to the above, we can see
$$
\lim{X_S}\Ext_{\text{graded }{\mathrm{gr}}^FD_S\text{-mod}}^p
(\O_S,{\mathrm{gr}}^FH_{dR}^q(X_S/S))\simeq
\Lambda_X^{r-p,q-r+p}(p),
$$
and the map \eqref{xtol} is compatible to the map between Yoneda
extension groups induced from the exact functor
${\mathrm{MF}}_{rh}(S)\ra \{\text{graded }{\mathrm{gr}}^FD_S\text{-mod}\}$,
$(M,F)\mapsto {\mathrm{gr}}^FM$.

\medskip

Together with Proposition \ref{xi527}, we thus have:
\begin{prop}\label{muminv58}
Let $X$ be a nonsingular projective variety over $\C$.
Then we have the natural map from the extension groups in $\MM(\C)$ to
the space of Mumford's infinitesimal invariants
as the composition of the map in Proposition $\ref{xi527}$ and \eqref{xtol}:
$$
\Ext_{\MM(\C)}^p(\Q(0),H^q(X)(r))\lra \Lambda_X^{r-p,q-r+p}(p).
$$
\end{prop}
\section{Algebraic cycles and extensions of arithmetic Hodge structure}
\label{anobservationtocycles}
In this section, we study the extension groups in the abelian category
$\MM(\C)$ and algebraic cycles, in particular, the Bloch conjecture.

\medskip

\subsection{Higher Abel-Jacobi maps}
Let $X$ be a nonsingular variety over $\C$.
We denote its higher Chow group by $\CH^r(X,m)$.
Then we have the following map
\begin{equation}\label{cyclemap58}
\begin{CD}
\lim{X_S}\CH^r(X_S,m)@>>> \lim{X_S}
\Ext^{2r-m}_{\MHM(X_S)}(\Q_{X_S}(0),\Q_{X_S}(r))\\
\Vert@.\Vert\\
\CH^r(X,m)@>>>
\Ext^{2r-m}_{\MM(X)}(\Q_{X}(0),\Q_{X}(r)),
\end{CD}
\end{equation}
where the top horizontal arrow is the one in Corollary \ref{cycle529}.
We denote the above map by $c=c^r(X,m)$.
On the other hand, there is the Leray spectral sequence
\begin{equation}\label{leraymc}
E_2^{pq}=\Ext^p_{\MM(\C)}(\Q(0),H^q(X)(r))\Longrightarrow
\Ext^{p+q}_{\MM(X)}(\Q_X(0),\Q_X(r)),
\end{equation}
as the inductive limit of the Leray spectral sequence for the model $X_S$
(Corollary \ref{leray58}).
It degenerates at $E_2$ terms if $X$ is projective.
The spectral sequence
\eqref{leraymc} defines the Leray filtration on the extension group
$\Ext^{\bullet}_{\MM(X)}(\Q_X(0),\Q_X(r))$, which we denote by $F^{\bullet}$.

\begin{definition}
Let $X$ be a nonsingular projective variety over $\C$.
We define the filtration on the higher Chow group as follows:
$$F^{\nu}\CH^r(X,m):=c^{-1}F^{\nu}
\Ext^{2r-m}_{\MM(X)}(\Q_{X}(0),\Q_{X}(r)).$$
Then, together with the $E_2$-degeneration of the spectral sequence 
\eqref{leraymc}, we obtain:
$$
\rho^{\nu}_X:\Gr^{\nu}_F\CH^r(X,m)\lra 
\Ext^{\nu}_{\MM(\C)}(\Q(0),H^{2r-m-\nu}(X)(r)).
$$
We call the above the {\it $\nu$-th higher Abel-Jacobi map}.
\end{definition}

By definition, 
we can easily show Theorem \ref{intro58} \eqref{fil531} 
and \eqref{filt531}. For example, Theorem \ref{intro58} \eqref{filt531} (d)
can be shown by 
using the argument in \cite{murre2}. 
Theorem \ref{intro58} \eqref{filt531} (e) follows from the hard Lefschetz
theorem in the terminology of arithmetic Hodge structure. 
(See \cite{A} for details.)
\begin{example}
Let $r=\dim X=n$, $\nu=2$ and $m=0$. Then the higher Abel-Jacobi map gives
the following map:
\begin{equation}\label{alb2}
\rho_X^2:T(X)\lra \Ext^{2}_{\MM(\C)}(\Q(0),H^{2n-2}(X)(n)).
\end{equation}
We call the above \eqref{alb2} the {\it second Albanese map}.
\end{example}
\subsection{Mumford's infinitesimal invariants}
Composing the higher Abel-Jacobi map and the maps in 
Proposition \ref{xi527} and \ref{muminv58},
we have Theorem \ref{intro58} \eqref{xiandde531}:
$$
\xi^{\nu}_X:\Gr^{\nu}_F\CH^r(X,m)
\lra \X_X^{r-\nu,r-m}(\nu),
$$
$$
\delta^{\nu}_X:\Gr^{\nu}_F\CH^r(X,m)
\lra \Lambda_X^{r-\nu,r-m}(\nu).
$$

When $m=0$, the above maps have already been defined in \cite{AS}.

\medskip

We know the map of infinitesimal invariants $\delta_X$
is not zero in general.
More strongly, the following holds.
\begin{theorem}[\cite{AS} Theorem  (0-6)]
Let $X\subset \P^{m+s}_{\C}$ be a nonsingular complete intersection variety
of degree $(d_1,\cdots,d_s)$, with a
sufficiently general coefficients. Assume that
one of the following conditions holds.
\begin{enumerate}
\renewcommand{\labelenumi}{(\theenumi)}
\item
$m\geq2$ and ${\mathrm{min}}\{d_i\}\cdot s\geq
2m+s.$
\item
$m=1$ and  ${\mathrm{min}}\{d_i\}\cdot s\geq 3+s.$
\end{enumerate}
Then the map of infinitesimal invariant induces an injective map
\begin{equation}\label{shujisaito61}
\delta_X^{m}:F^{m}\CH^m(X)_{\mathrm{gen}}\lra \Lambda^{0,m}_{X}(m).
\end{equation}
Here $F^{m}\CH^m(X)_{\mathrm{gen}}$ denotes the subgroup of $F^{m}\CH^m(X)$ 
generated
by $0$-cycles in general positions $($loc.cit. Ch.II Definition 
$(5$-$4$-$1))$.
In particular, we have 
$$
F^{m+1}\CH^m(X)_{\mathrm{gen}}=0.
$$
\label{shujisaito}
\end{theorem}
\begin{remark}
For a nonsingular complete intersection $X$ of dimension $m$, it is easy to
see that
$F^2\CH^m(X)=\cdots=F^m\CH^m(X)=T(X)$ and
$\delta_X^{m}(F^{m+1}\CH^m(X))=0$, where $F^{\bullet}$ denotes either
our filtration or S.Saito's filtration
(\cite{AS} Ch.II \S2, \cite{shuji}, \cite{ss}).
It is conjectured that our filtration coincides with S.Saito's filtration.
But I do not have the proof of it.
\end{remark}
\begin{pf}(Sketch).
We prove the assertion by the induction of $\dim X$.

We first show the map
\begin{equation}\label{monod1}
\delta_C^{1}:F^{1}\CH^1(C)\lra \Lambda^{0,1}_{C}(1).
\end{equation}
is injective for a general complete intersection curve 
$C\subset\P^{s+1}_{\C}$ of degree $(d_1,\cdots,d_s)$, which
satisfies ${\mathrm{min}}\{d_i\}\cdot s\geq 3+s$.
Let $f:C_S\ra S$ be a model of $C$, where we may assume that 
there is a dominant morphism $S\ra {\mathcal M}$ to the moduli space of
complete intersection curves of degree $(d_1,\cdots,d_s)$, because
the coefficients of $C$ is general.
Let $\Lambda^{0,1}_{C_S/S}(1):=\vg(S,
R^1f_*\O_{C_S}\ot 
\Omega_{S/\Qb}^1/\ol{\nb}(f_*\Omega^1_{C_S/S}))$.
Then the map \eqref{monod1} can be obtained by the inductive limit over $S$
of the map $F^{1}_S\CH^{1}(C_S)\ra \Lambda^{0,1}_{C_S/S}(1)$.
Here $F^{\bullet}_S\CH^{\bullet}(C_S)$ denotes the relative filtration
on the Chow group (cf. \cite{AS} Ch.II \S2).
Let $\Lambda^{0,1}_{C_S^{an}/S^{an}}(1):=\Lambda^{0,1}_{C_S/S}(1)
\ot\O_{S}^{an}\simeq\vg(S^{an},
R^1f_*\O_{C_S}\ot 
\Omega_{S^{an}/\C}^1/\ol{\nb}(f_*\Omega^1_{C_S^{an}/S^{an}}))$.
In order to prove the injectivity of \eqref{monod1}, it suffices to
show that if a cycle $z\in F^1_S\CH^1(C_S)$ has a zero-image under
the following map
\begin{equation}\label{monod2}
\delta_S^1:
F^1_S\CH^1(C_S)\lra \Lambda^{0,1}_{C_S^{an}/S^{an}}(1),
\end{equation}
then $z\vert _C$ is rationally equivalent to zero. 
Since $f_*\Omega^1_{C_S^{an}/S}\ot \Omega_{S^{an}/\C}^1\ra
R^1f_*\O_{C_S}\ot \Omega_{S^{an}/\C}^2$ is injective by
the symmetrizer lemma for complete intersections 
(loc.cit.Ch.I.Corollary (4-5)), we have the injection
\begin{align*}
\X_{C_S^{an}/S^{an}}^{0,1}(1)&:=\vg(S^{an},
(R^1f_*\Omega^{\bullet}_{C_S^{an}/S^{an}})\ot 
\Omega_{S^{an}/\C}^1/{\nb}(f_*\Omega^1_{C_S^{an}/S^{an}}))\\
&\hra \Lambda^{0,1}_{C_S^{an}/S^{an}}(1).
\end{align*}
There is the following exact sequence
$$
0\ra
\vg(S^{an},R^1f_*\C_{C_S^{an}}/R^1f_*\Q_{C_S^{an}})
\ra 
\vg(S^{an},J(C_S/S))
\ra
\X_{C_S^{an}/S^{an}}^{0,1}(1),
$$
where we put $J(C_S/S):=
R^1f_*\Omega^{\bullet}_{C_S^{an}/S^{an}}/
f_*\Omega^1_{C_S^{an}/S^{an}}+R^1f_*\Q_{C_S^{an}}
$.
It is easy to see that 
$\vg(S^{an},R^1f_*\C_{C_S^{an}})=0$ and hence 
$\vg(S^{an},R^1f_*\C_{C_S^{an}}/R^1f_*\Q_{C_S^{an}})=0$ because of
the dominant morphism $S\ra {\mathcal M}$.
Thus it follows from the Abel's theorem that if $\delta_S(z)=0$, then
$z\vert_{f^{-1}(t)}=0$ for any $t\in S(\C)$, in particular,
$z\vert_{C}=0$.

\medskip

Next we prove for any dimensional case. Let $X$ be
as in Theorem \ref{shujisaito} and $m\geq2$.
 let 
$$X=X_m\supset X_{m-1}\supset\cdots\supset X_2\supset X_1=C$$ 
be a sequence of general hyperplane sections, in which the multi-degree
of each $X_i$ satisfies the condition in Theorem \ref{shujisaito}.
Then there is the commutative diagram
$$
\begin{CD}
F^{m}\CH^m(X)@>{\delta_X^m}>> \Lambda^{0,m}_{X}(m)\\
@AAA@AA{\alpha_{m-1}}A\\
F^{m-1}\CH^{m-1}(X_{m-1})@>{\delta_{X_{m-1}}^{m-1}}>> 
\Lambda^{0,m-1}_{X_{m-1}}(m-1)\\
@AAA@AA{\alpha_{m-2}}A\\
\vdots @. \vdots\\
@AAA@AA{\alpha_{1}}A\\
F^{1}\CH^{1}(C)@>{\delta_C^{1}}>> 
\Lambda^{0,1}_{C}(1).
\end{CD}
$$
By the above, $\delta_C^1$ is injective. 
Moreover, 
each $\alpha_i$ is also injective
by the symmetrizer lemma for open complements
(loc.cit. Ch.I Theorem (4-4)).
Put 
$$F^{m}\CH^m(X)_{\mathrm{gen}}:=\Image(\sum_CF^{1}\CH^{1}(C)\ra 
F^{m}\CH^{m}(X))$$ where $C$ runs over all curves as above.
Then we have that the map \eqref{shujisaito61}
$${\delta_X^m}:
F^{m}\CH^m(X)_{\mathrm{gen}}\lra \Lambda^{0,m}_{X}(m)
$$
is injective.
\end{pf}
Mumford's infinitesimal invariants have the advantage
of explicit calculations.
But it is a ``coarse" invariant, that is,
there are several examples of an algebraic cycle which is not rationally 
equivalent
to zero, 
but its infinitesimal invariant is zero (cf. \cite{G2}).
We hope that our higher Abel-Jacobi invariants capture all cycles
(cf. Conjecture \ref{ext2} below).
In fact,
there is an example of a $0$-cycle on a surface such that 
its infinitesimal invariant vanishes, but 
the higher Abel-Jacobi invariant does not:

\begin{theorem}\label{asa68}
Let $C$ be a projective nonsingular curve over $\Qb$ 
with genus $g\geq2$, such 
that $\rank \NS(C\times C)=3$
\footnote{
The proof of the existence of such a curve was taught from S.Mochizuki
(the private communication with the author).}.
Let $O\in C(\Qb)$ be a $\Qb$-valued point such that $K_C-(2g-2)\cdot O$
is not $\Q$-linearly equivalent to $0$
\footnote{
There are infinitely many such 
points $O$ by Raynaud's theorem (\cite{Ray}).}.
Let $P\in C(\C)\setminus C(\Qb)$ be any $\C$-valued point which is not
$\Qb$-valued one.
Put $X:=C_{\C}\times C_{\C}$ and $z:=(P,P)-(P,O)-(O,P)+(O,O)\in T(X)$.
Then we have
$$
\xi_X^{2}(z)=0 \quad \text{in }\quad
\X^{0,2}_X(2),
$$
but
$$
\rho_X^2(z)\not=0 \quad \text{in }
\Ext_{\MM(\C)}^2(\Q(0),H^2(X)(2)).
$$
\end{theorem}
The proof of the above theorem will appear in the forthcoming
paper \cite{A}.

\subsection{Bloch conjecture}
Recall the famous {\it Bloch conjecture} (\cite{bloch} Lec.1):
\begin{claim}[Bloch]\label{Bloch}
Let $X$ be a projective nonsingular surface over $\C$ with the geometric
genus $p_g=\dim H^0(K_{X/\C})=0$. Then the $($usual$)$ Albanese map
$$
A_0(X) \lra \Alb(X)
$$
is injective.
\end{claim}
This conjecture was observed by S.Bloch with an ingenious insight, and
it led to the Bloch-Beilinson formula~\eqref{blochbeilinson}.
The Bloch conjecture was proved for surfaces which are not of general type
(\cite{bkl}).
But, in general, it is still open.
\vspace{5mm}

Now we conjecture:
\begin{claim}\label{ext2}
The second Albanese map~\eqref{alb2}
$$\rho_X^2:
T(X) \lra \Ext^2_{\MM(\C)}(\Q(0),H^{2}(X)(2))
$$
is injective for any projective nonsingular
surface $X$ $($not necessarily $p_g=0)$.
\end{claim}
Or, more generally
\begin{claim}\label{extn}
The cycle map $c=c^r(X,m)$ in \eqref{cyclemap58} is injective
for any projective nonsingular variety $X$ over $\C$.
\end{claim}

\medskip

Conjecture \ref{ext2} implies the 
Bloch conjecture (use the method in \cite{bloch} Lec. 1).
Moreover Conjecture \ref{extn} 
also has an application to algebraic cycles on varieties over
an algebraic number field.
Let us explain it.

Let $X_{\Qb}$ be a nonsingular projective variety over $\Qb$.
By the construction, the higher Abel-Jacobi map factors as follows:
$$
\Gr_F^{\nu}\CH^r(X_{\Qb}) \ra \Ext^{\nu}_{\MHM(\Qb)}(\Q(0),H^{2r-\nu}(X)(r))
\ra \Ext^{\nu}_{\MM(\C)}(\Q(0),H^{2r-\nu}(X)(r)).
$$
Here $F^{\nu}\CH^r(X_{\Qb})$ denotes the subgroup of 
$F^{\nu}\CH^r(X_{\C})$ generated by
cycles of codimension $r$ defined over $\Qb$.
Since the arithmetic Gauss-Manin connection in $\MHM(\Qb)$
vanishes, we can show any extension groups in $\MHM(\Qb)$ 
vanishes if $\nu\geq2$ (cf. Theorem~\ref{carl} (2)).
Then, due to Conjecture~\ref{extn}, we have
$$
F^{2}\CH^r(X_{\Qb}) =0.
$$
In particular, 
$$
T(X_{\Qb})=0.
$$
Let $X_k$ be a nonsingular projective variety over an algebraic
number field $k$.
Due to the injectivity of the map $\CH_0(X_k)\ra\CH_0(X_{\Qb})$ (note
that the Chow groups are always tensored with $\Q$) , 
we also have $T(X_k)=0$. Combined with the Mordell-Weil theorem,
we obtain the finiteness of the rank of the Chow group:
$$
{\mathrm{rank}}~\CH_0(X_k)< \infty.
$$

\begin{theorem}
\begin{enumerate}
\renewcommand{\labelenumi}{(\theenumi)}
\item\label{candq1531}
Conjecture $\ref{ext2}$ implies the Bloch conjecture $\ref{Bloch}$.
\item\label{candq2531}
Let $X_k$ be a nonsingular projective variety over an algebraic
number field $k$. Assume that Conjecture $\ref{extn}$ holds.
Then we have
$ F^{2}\CH^r(X_{\Qb}) =0$.
In particular, it implies that the rank of $\CH_0(X_k)$ is finite.
\end{enumerate}
\label{candq531}
\end{theorem}

Further we can show the `almost' 
converse of Theorem \ref{candq531} \eqref{candq2531} for $r=2$:
\begin{theorem}[\cite{A}]
Assume the followings hold:
\begin{enumerate}
\renewcommand{\labelenumi}{(\theenumi)}
\item\label{bei1531}
The Abel-Jacobi map $($cf. \S$2.5)$
$$
\rho:\CH^2(X_{\Qb})_{\hom}\lra J^2(X_{\C})
$$
is injective for any projective nonsingular
variety $X_{\Qb}$ of arbitrary dimension over $\Qb$.
$($Here we put $X_{\C}=X_{\Qb}\ot\C$.$)$
\item\label{bei2531}
Let $X_{\Qb}$ and $\rho$ be as above.
Let $z\in \CH^2(X_{\C})_{\hom}$ be any algebraic cycle.
Then $\rho(z)=0$ if and only if $\rho(z^{\sigma})=0$
for any $\sigma\in \Aut(\C/\Qb)$.
$($Note that $\Aut(\C/\Qb)$ acts on $\CH^2(X_{\C})_{\hom}$.$)$
\end{enumerate}
Then Conjecture $\ref{extn}$ for $r=2$ and $m=0$ is true. In particular,
\eqref{bei1531} and \eqref{bei2531} implies the Bloch conjecture 
$\ref{Bloch}$.
\label{bei531}
\end{theorem}

Theorem \ref{bei531} \eqref{bei1531} is nothing
but the conjecture due to
A.Beilinson.
Theorem \ref{bei531} \eqref{bei2531} should be true
 if $F^2_{\mathcal M}=\ker \rho$.
 
\subsection{}
Finally we propose the following 
problem about surjectivity of the cycle map.

Let $X$ be a projective nonsingular variety over $\C$, and put
$$J^r(X)_{\mathrm{arith}}=\Image(\Ext^1_{\MM(\C)}(\Q(0),H^{2r-1}(X)(r))
\lra J^r(X)_{\Q})$$ for $r\geq2$ (cf. \S\ref{carlson}).
\begin{exmp}
Is the Abel-Jacobi map 
$$\CH^r(X)_{\mathrm{hom}}\lra J^r(X)_{\mathrm{arith}}$$
surjective?
\end{exmp}
In general, the Abel-Jacobi map \eqref{AJ} is not surjective.
For example, when the Hodge number 
$\Sigma_{k\geq r+1}h^{k,2r-1-k}(X)\not=0$, 
it is not so
(or, more strong result, see \cite{green}). 
It is an interesting problem to give Hodge
theoretic characterization of its image.

\end{document}